\documentclass[12pt,leqno]{article}
\usepackage{amsfonts, color,amsmath}

\newcounter{conjecture}
\newcounter{remark}

\newtheorem{theorem}{Theorem}
\newtheorem{lemma}{Lemma}
\newtheorem{proposition}{Proposition}

\newcommand{\lar}{\longrightarrow}
\newcommand{\eps}{\varepsilon}

\newcommand{\nats}{\mathbb{N}}
\newcommand{\lll}{\label}
\newcommand {\rrr}[1]{(\ref{#1})}

\def \be{\begin{equation}}
\def \ee{\end{equation}}
\def \bt{\begin{theorem}}
\def \et{\end{theorem}}
\def \bc{\begin{corollary}}
\def \ec{\end{corollary}}
\def \bl{\begin{lemma}}
\def \el{\end{lemma}}
\def \bea{\begin{eqnarray}}
\def \eea{\end{eqnarray}}
\def \bas{\begin{eqnarray*}}
\def \eas{\end{eqnarray*}}

\def \ga{\gamma}

\def \vski{\vspace{12pt}}
\def \ff{\infty}

\def \DD{\Delta}

\def \({\left(}
\def \){\right)}

\def \nn{\nonumber}

\def \bc{\begin{center} }
\def \ec{\end{center} }
\def \bs{\begin{slide} }
\def \es{\end{slide} }

\def\square{{\vcenter{\vbox{\hrule height.3pt
         \hbox{\vrule width.3pt height5pt \kern5pt
            \vrule width.3pt}
         \hrule height.3pt}}}}
\def\qed{{\hfill $\square$ \bigskip}}

\newcounter{cccases}
\begin{document}

\title{An Optional Stopping-like theorem for large stopping times on birth-death chains.}

\author{
\begin{tabular}{c}
\textit{Greg Markowsky} \\
gmarkowsky@gmail.com \\
+61 03 9905-4487 \\
Monash University \\
Department of Mathematical Sciences \\
Victoria, 3800 Australia
\end{tabular}}


\maketitle

\begin{abstract}
An elementary proof is given for a theorem showing that certain birth-death chains show martingale-like behavior at large stopping times. This is a generalization of and new proof for a theorem from \cite{mebaby}.
\end{abstract}

\section{Introduction}

Let $X_m$ be a Markov chain taking values on the nonnegative integers with the following transition probabilities for $n \neq 0$

\be \lll{castle} p_{nj} = \left \{ \begin{array}{ll}
r_{n} & \qquad  \mbox{if } j=n+1  \\
l_n & \qquad \mbox{if } j=n-1 \\ 0 & \qquad \mbox{if } |n-j| \neq 1\;.
\end{array} \right. \ee

Implicit here is the fact that $r_n+l_n=1$. We suppose further that $X_0 = k$ almost surely, for some $k \in \nats$, so that at time $0$ our chain is at a fixed state. To avoid trivialities, assume $l_n, r_n>0$ for all $i$. $X_m$ is essentially a random walk on the nonnegative integers, moving to the right from state $n$ with probability $r_n$ and to the left with probability $l_n$. Such a Markov chain is referred to as a {\it birth-death chain}\footnotemark. This name comes from considering $X_m$ as the number of members in a population, where at each step either a new member is born or an old member dies, causing the process to increase or decrease by 1. We can assume $p_{00}=1$ and $p_{0j}=0$ for any $j \neq 0$, as when the population reaches 0 it is considered to have gone extinct with no possibility of regeneration. We let $T_\DD$ denote the first time $X_m$ hits $0$, with $T_\DD=\ff$ if it never does. If $l_n = r_n = 1/2$ for all $n$, then our birth-death chain is simple random walk stopped at 0, and is therefore a martingale. The well-known Optional Stopping Theorem then implies that $E[X_\tau]= E[X_0]=k$ for any stopping time $\tau$ such that $E[\tau]< \ff$.

\footnotetext{One may also find the name {\it state dependent random walk} in the literature.}

\vski

Clearly, a general birth-death chain is not a martingale, and therefore $E[X_\tau]= E[X_0]$ need not hold. However, we can give a sufficient condition such that similar behavior holds for large stopping times. The following is the main result of this paper.

\bt \lll{bigguy}
If $t_\ff := \lim_{n \lar \ff} \frac{l_1\ldots l_{n}}{r_1\ldots r_{n}}$ exists then $\lim_{m \lar \ff}E[X_{T_m}]$ exists, and

\be \lll{deb}
\lim_{m \lar \ff}E[X_{T_m}] = \frac{1+ \frac{l_1}{r_1} + \ldots + \frac{l_1\ldots l_{k-1}}{r_1\ldots r_{k-1}}}{t_\ff}
\ee

for any non-decreasing sequence of stopping times $\{T_m\}_{m=1}^\ff$ satisfying $E[T_m]<\ff$ such that $T_m \lar T_\DD$ almost surely.
\et

Note that the theorem is valid if $t_\ff = +\ff$ or $0$, with \rrr{deb} being interpreted as $0$ or $\ff$ respectively. The condition $E[T_m]<\ff$ is clearly necessary, since in the case $P(T_\DD < \ff)=1$ we have $E[X_{T_\DD}]=0$ regardless of the value of $t_\ff$. A preliminary version of this theorem has already been proved in \cite{mebaby}; only the case $T_m = m \wedge T_\DD$ was considered, but the extension to arbitrary stopping times causes no real difficulty, as will be shown below. The proof given in \cite{mebaby} relies on properties of Brownian
motion and its corresponding theory of local time. Given the simplicity of the statement of the theorem, it seems desirable to have a proof which is elementary and self-contained, without reference to the machinery of It\={o} calculus and the like. The purpose of this paper is to give such a proof, as well as the extension to more general stopping times than $T_m = m \wedge T_\DD$.

\section{Proof of Theorem}

If we condition on $X_{m'-1}=n$, then if we move to the right we increase the expectation by 1, and if we move to the left we decrease by 1. It follows that

\be \label{}
E[X_{m'}]-E[X_{m'-1}] = \sum_{n=1}^{\ff} P(X_{m-1}=n)(r_n-l_n)
\ee

Summing over all $0 \leq m' \leq m$ gives

\be \label{nuts}
E[X_{m}]-k=E[X_{m}]-E[X_0] = \sum_{n=1}^{\ff} E[G_{m}^n](r_n-l_n)
\ee

where $G_{m}^n$ is the number of times that $X$ is equal to $n$ on or before time $m-1$. We need the analogous formula for stopping times.

\begin{lemma} \label{kris}
Let $T$ be a stopping time with $E[T]<\ff$. Then

\begin{equation} \label{}
E[X_{T}] = k + \sum_{n=1}^{\ff} E[G_{T}^n](r_n-l_n)
\end{equation}

where $G_{T}^n$ is the number of times that $X$ is equal to $n$ on or before time $T-1$.
\end{lemma}

{\bf Proof:} For any integer $h \geq 0$ let $S_h = T \wedge h$. For $h' \geq 1$ we have

\begin{equation} \label{}
E[X_{S_{h'}} - X_{S_{h'-1}}] = \sum_{n=1}^{\ff} P(T>h'-1,X_{h'-1}=n)(r_n-l_n)
\end{equation}

Sum this expression over $1 \leq h' \leq h$ and use $E[X_{S_0}]=k$ to get

\begin{equation} \label{}
E[X_{S_{h}}] = k + \sum_{h'=1}^{h} \sum_{n=1}^{\ff} P(T>h'-1,X_{h'-1}=n)(r_n-l_n)
\end{equation}

As $S_h \leq h$ we have $\sup_{0 \leq m \leq S_h}X_m \leq k+h$, and $\{T>h'-1\} = \{S_{h} > h'-1\}$. Thus, the following manipulations are justified.

\begin{equation} \label{xinyi}
\begin{split}
E[X_{S_{h}}] & = k + \sum_{h'=1}^{h} \sum_{n=1}^{k+h} P(S_h>h'-1,X_{h'-1}=n)(r_n-l_n) \\
& = k + \sum_{n=1}^{k+h} (r_n-l_n) \sum_{h'=1}^{h} P(S_h>h'-1,X_{h'-1}=n) \\
& = k + \sum_{n=1}^{k+h} (r_n-l_n) E[ \sum_{h'=1}^{h} 1(\{S_h>h'-1,X_{h'-1}=n\})] \\
& = k + \sum_{n=1}^{k+h} E[G_{S_h}^n](r_n-l_n)
\end{split}
\end{equation}

where $1(A)$ denotes the indicator of the set $A$. Since $E[T] < \ff$, $\lim_{h \lar \ff}E[T-S_h] = 0$, and since $|X_{m+1}-X_m|=1$ almost surely we have $E[|X_T - X_{S_h}|] \leq E[T-S_h] \lar 0$ as well. The result therefore follows by letting $h \lar \ff$ in \rrr{xinyi}. \qed

Let $t_0 := 1$ and

\be \label{pred}
t_n := \frac{l_1 l_2 \ldots l_n}{r_1 r_2 \ldots r_n}
\ee

for $n>0$. The following lemma gives the exit probabilities of an interval.

\begin{lemma} \label{sound}
Suppose $Y_m$ is a birth-death chain with the same transition probabilities as $X_m$ but with $Y_0=k'$ almost surely. Let $0 \leq a < k' < b < \ff$, and let $\tau$ be the first time that $Y_m$ hits $a$ or $b$. Then

\be \label{}
P(Y_\tau=a) = \frac{\sum_{n=k'}^{b-1}t_n}{\sum_{n=a}^{b-1}t_n}, \; \; P(Y_\tau=b) = \frac{\sum_{n=a}^{k'-1}t_n}{\sum_{n=a}^{b-1}t_n}
\ee
\end{lemma}

There are a number of ways to prove this well-known fact. One way is to write down the correct recurrence relations and verify that the given expressions satisfy them; see \cite{norris} or \cite{syski} for the details. Another is to use the relationship between electrical resistance on graphs and random walks; see \cite{doysne} for an elegant introduction to this technique. Yet a third is to realize the birth-death chain as a Brownian motion evaluated at a properly chosen sequence of increasing stopping times; see \cite{mebaby}. Setting $a=0$ and letting $b \lar \ff$ we obtain

\be \label{seat}
P(X_m = 0 \mbox{ for some }m) = \frac{\sum_{n=k}^\ff t_n}{\sum_{n=0}^\ff t_n}
\ee

with the quotient interpreted as equal to $1$ if the sums diverge. Note that if $P(X_m = 0 \mbox{ for some }m) < 1$ then the sums in \rrr{seat} converge so that $t_n \lar 0$ as $n \lar \ff$. Furthermore, in that case we must have $X_{T_m} \lar \ff$ on some set of positive measure, so that $E[X_{T_m}] \lar \ff$, and we see that \rrr{deb} holds trivially. We will therefore assume from now on that $P(X_m = 0 \mbox{ for some }m) = 1$. In what follows, let $x_n = \sum_{i=0}^{n-1} t_i$. Lemma \ref{kris} shows the importance of calculating $E[G_{T_\DD}^n]$, and we therefore need the following.

\begin{lemma} \label{beck}
\be \label{}
E[G_{T_\DD}^n] = \lim_{m \lar \ff} E[G_{T_m}^n] = \frac{\min(x_n,x_k)}{t_{n-1}l_n} 
\ee
\end{lemma}

{\bf Proof:} That $E[G_{T_\DD}^n] = \lim_{m \lar \ff} E[G_{T_m}^n]$ is clear by monotone convergence. Let $\ga_n = E[G_{T_\DD}^n]$, and let us first suppose $n=k$. If $X_1 = k+1$, then $X$ must eventually return to $k$, whereas if $X_1=k-1$ then by Lemma \ref{sound} the probability of returning to $k$ before hitting $0$ is $\frac{x_{k-1}}{x_k}$. Thus, $\ga_k$ satisfies

\be \label{}
\ga_k = 1 + r_k \ga_k + l_k \Big(\frac{x_{k-1}}{x_k}\Big) \ga_k
\ee

Solving for $\ga_k$ gives

\begin{equation*} \label{}
\ga_k = \frac{1}{1-r_k-l_k\Big(\frac{x_{k-1}}{x_k}\Big)} = \frac{x_k}{x_k(1-r_k)-x_{k-1}l_k} =\frac{x_k}{(x_k-x_{k-1})l_k}=\frac{x_k}{t_{k-1}l_k}
\end{equation*}

For general $n$, we let $E_j[G_{T_\DD}^n]$ denote the expected number of times that the birth-death chain $Y_m$ hits $n$ before hitting $0$, where $Y_m$ has the same transition probabilities as $X_m$ but $Y_0=j$ almost surely. Then by the strong Markov property of $X$ we have

\be
E[G_{T_\DD}^n] = P(X_m \mbox{ hits $n$ before 0})E_n [G_{T_\DD}^n] = P(X_m \mbox{ hits $n$ before 0})\frac{x_n}{t_{n-1}l_n}
\ee

The general result now follows by noting that $P(X_m \mbox{ hits $n$ before $0$})$ is 1 if $n \leq k$ and $\frac{x_k}{x_n}$ if $n>k$ by Lemma \ref{sound}. \qed

We see initially from this that $E[G_{T_m}^n] \leq \frac{x_k}{t_{n-1}l_n}$. Returning to \rrr{nuts}, we may write $(r_n-l_n) = t_{n-1}l_n(\frac{1}{t_n}-\frac{1}{t_{n-1}})$. If we assume $\sum_{n=1}^{\ff}|\frac{1}{t_n}-\frac{1}{t_{n-1}}| < \ff$ then by \rrr{nuts} we have

\be \label{}
\begin{split}
E[X_{T_m}] &= k + \sum_{n=1}^{\ff} E[G_{T_m}^n](r_n-l_n) \leq k + \sum_{n=1}^{\ff} \frac{x_k}{t_{n-1}l_n}t_{n-1}l_n \Big|\frac{1}{t_n}-\frac{1}{t_{n-1}}\Big| \\ & \leq x_k \sum_{n=1}^{\ff} \Big|\frac{1}{t_n}-\frac{1}{t_{n-1}}\Big| < \ff
\end{split}
\ee

We may then apply the dominated convergence theorem(to the functions $f_m(n) = E[G_{T_m}^n](r_n-l_n)$ defined on $\nats$) to conclude that $\lim_{m \lar \ff} E[X_{T_m}]$ exists and is given by

\be \label{}
\begin{split}
\lim_{m \lar \ff}E[X_{T_m}] & = k + \sum_{n=1}^{\ff} E[G_{T_\DD}^n](r_n-l_n) \\ & = k + \sum_{n=1}^{k} \frac{x_n}{t_{n-1}l_n}(r_n-l_n) + \sum_{n=k+1}^{\ff}\frac{x_k}{t_{n-1}l_n}(r_n-l_n) \\ & = k + \sum_{n=1}^{k} x_n(\frac{1}{t_n}-\frac{1}{t_{n-1}}) + x_k \sum_{n=k+1}^{\ff}(\frac{1}{t_n}-\frac{1}{t_{n-1}}) \\ & = k + x_k \frac{1}{t_k} - x_1\frac{1}{t_0} - \sum_{n=1}^{k-1} \frac{1}{t_n}(x_n - x_{n-1}) + x_k(\frac{1}{t_\ff} - \frac{1}{t_k})
\end{split}
\ee

where summation by parts(see \cite{lang}) was used in passing to the last line. It was also assumed that $0<t_\ff<\ff$, but the cases $t_\ff = 0, \ff$ are handled easily by the same technique. As $\frac{1}{t_n}(x_n - x_{n-1}) = x_1 \frac{1}{t_0} = 1$, this simplifies to

\be \label{}
\lim_{m \lar \ff}E[X_{T_m}] = \frac{x_k}{t_\ff} = \frac{1+ \frac{l_1}{r_1} + \ldots + \frac{l_1\ldots l_{k-1}}{r_1\ldots r_{k-1}}}{t_\ff}
\ee

which is the desired expression. Now we must remove the assumption $\sum_{n=1}^{\ff}|\frac{1}{t_n}-\frac{1}{t_{n-1}}| < \ff$. We need the following lemma.

\begin{lemma} \label{ddd}
If $\tau$ is any stopping time, then $E[G_\tau^n]t_{n-1}l_n > E[G_\tau^{n+1}]t_{n}l_{n+1}$ for any $n \geq k$.
\end{lemma}

{\bf Proof:} Using the strong Markov property of $X_m$ and Lemma \ref{beck}, we have

\begin{equation} \lll{cas1}
\begin{split}
E[G_\tau^n]t_{n-1}l_n & = E[G_{T_\DD}^n]t_{n-1}l_n - E[\sum_{m=\tau+1}^{T_\DD} 1_{\{n\}}(X_m)]t_{n-1}l_n \\
& = x_k - t_{n-1}l_n\sum_{j=1}^{\ff} P(X_{\tau+1}=j) E_j[G_{T_\DD}^n] \\
& = x_k - t_{n-1}l_n \sum_{j=1}^{n} P(X_{\tau+1}=j) \frac{x_j}{t_{n-1}l_n} - t_{n-1}l_n\sum_{j=n+1}^{\ff} P(X_{\tau+1}=j) \frac{x_n}{t_{n-1}l_n} \\
& = x_k - \sum_{j=1}^{n} P(X_{\tau+1}=j) x_j - \sum_{j=n+1}^{\ff} P(X_{\tau+1}=j) x_n
\end{split}
\end{equation}

Likewise,

\begin{equation} \label{cas2}
E[G_\tau^{n+1}]t_{n}l_{n+1} = x_k - \sum_{j=1}^{n+1} P(X_{\tau+1}=j) x_j - \sum_{j=n+2}^{\ff} P(X_{\tau+1}=j) x_{n+1}
\end{equation}

The result now follows from \rrr{cas1}, \rrr{cas2}, and the fact that $x_{n+1}>x_n$. \qed

Let us suppose initially that for each $m$ we can find a number $N_m>k$ such that $\sup_{0 \leq m' \leq T_m} X_{m'} \leq N_m$ almost surely. Returning to \rrr{nuts}, we have

\begin{equation} \label{rre}
E[X_{T_m}] = k + \sum_{n=1}^{k} E[G_{T_m}^n]t_{n-1}l_n(\frac{1}{t_n}-\frac{1}{t_{n-1}}) + \sum_{n=k+1}^{N_m+1} E[G_{T_m}^n]t_{n-1}l_n(\frac{1}{t_n}-\frac{1}{t_{n-1}})
\end{equation}

since $E[G_{T_m}^n] = 0$ for $n>N_m$. 
As $E[G_{T_m}^n] \lar E[G_{T_\DD}^n]$, the first sum converges to

\be \lll{true}
\begin{split}
\sum_{n=1}^{k} E[G_{T_\DD}^n]t_{n-1}l_n(\frac{1}{t_n}-\frac{1}{t_{n-1}}) & = \sum_{n=1}^{k} x_n (\frac{1}{t_n}-\frac{1}{t_{n-1}}) \\ & = x_k \frac{1}{t_k} - x_1 \frac{1}{t_0} - \sum_{n=1}^{k-1} \frac{1}{t_n}(x_{n+1}-x_n) \\ & = x_k \frac{1}{t_k} - k
\end{split}
\ee

We therefore need only show that the second sum converges to the proper thing. Let us assume that $t_\ff \in (0, \ff)$, and let $\eps>0$ be given. By summation by parts, using $E[G_{T_m}^{N_m+1}] = 0$, we have

\bea \label{dofo}
&& \sum_{n=k+1}^{N_m+1} E[G_{T_m}^n]t_{n-1}l_n(\frac{1}{t_n}-\frac{1}{t_{n-1}}) = 
\\ \nn && \hspace{1cm} - \frac{1}{t_k} E[G_{T_m}^{k+1}]t_{k}l_{k+1} + \sum_{n=k+1}^{N_m}\frac{1}{t_n}( E[G_{T_m}^n]t_{n-1}l_n - E[G_{T_m}^{n+1}]t_{n}l_{n+1})
\eea

Note that Lemma \ref{ddd} implies that the terms in the sum on the right side of \rrr{dofo} are all positive. We may choose $N>k$ such that $\frac{1}{t_n} \in (\frac{1}{t_\ff}-\eps,\frac{1}{t_\ff}+\eps)$ for all $n \geq N$. Having chosen this, we may choose $M > N-k$ such that $x_k \geq E[G_{T_m}^n]t_{n-1}l_n > x_k - \eps$ for all $n \in [k,N+1], m \geq M$. Setting $\frac{1}{t}^* = \sup_{j > 0} \frac{1}{t_j}$ and using \rrr{dofo}, we see that for $m>M$

\bea \label{}
&& \sum_{n=k+1}^{N_m} E[G_{T_m}^n]t_{n-1}l_n(\frac{1}{t_n}-\frac{1}{t_{n-1}})
\\ \nn && \hspace{1cm} \leq -\frac{1}{t_k} (x_k - \eps) + \sum_{n=k}^{N} \frac{1}{t}^* (E[G_{T_m}^n]t_{n-1}l_n - E[G_{T_m}^{n+1}]t_{n}l_{n+1})
\\ \nn && \hspace{2cm} + \sum_{n=N+1}^{N_m} (\frac{1}{t_\ff}+\eps) (E[G_{T_m}^n]t_{n-1}l_n - E[G_{T_m}^{n+1}]t_{n}l_{n+1})
\\ \nn && \hspace{1cm} \leq -\frac{1}{t_k} (x_k - \eps) + \frac{1}{t}^* (E[G_{T_m}^{k+1}]t_{k}l_{k+1} - E[G_{T_m}^{N+1}]t_{N}l_{N+1})
\\ \nn && \hspace{2cm} + (\frac{1}{t_\ff} + \eps) (E[G_{T_m}^{N+1}]t_{N}l_{N+1} - E[G_{T_m}^{N_m+1}]t_{N_m}l_{N_m+1})
\\ \nn && \hspace{1cm} \leq -\frac{1}{t_k} (x_k - \eps) + \frac{1}{t}^* \eps + (\frac{1}{t_\ff} + \eps) x_k \Big)
\eea

As $\eps>0$ is arbitrary, this shows that

\be \lll{}
\limsup_{m \lar \ff} \sum_{n=k+1}^{N_m} E[G_{T_m}^n]t_{n-1}l_n(\frac{1}{t_n}-\frac{1}{t_{n-1}}) \leq x_k(\frac{1}{t_\ff} - \frac{1}{t_k})
\ee

Proceeding similarly, we can obtain

\be \lll{}
\liminf_{m \lar \ff} \sum_{n=k+1}^{N_m} E[G_{T_m}^n]t_{n-1}l_n(\frac{1}{t_n}-\frac{1}{t_{n-1}}) \geq x_k(\frac{1}{t_\ff} - \frac{1}{t_k})
\ee

We conclude that

\be \lll{}
\lim_{m \lar \ff} \sum_{n=k+1}^{N_m} E[G_{T_m}^n]t_{n-1}l_n(\frac{1}{t_n}-\frac{1}{t_{n-1}}) = x_k(\frac{1}{t_\ff} - \frac{1}{t_k})
\ee

Combining this with \rrr{true} and \rrr{rre} gives the desired result. We assumed $\frac{1}{t_\ff} \in (0,\ff)$, but the cases $\frac{1}{t_\ff} \in \{0,\ff\}$ are similar but easier and are omitted. At this point we have but a single assumption which remains to be removed, which is that for each $m$ we can find a number $N_m>k$ such that $\sup_{0 \leq m' \leq T_m} X_{m'} \leq N_m$ almost surely. Suppose now that $T_m$ is a sequence of stopping times as in the statement of the theorem. As $E[T_m]<\ff$, there is a number $H_m$ such that if $S_m = T_m \wedge H_m$, then $E[T_m -S_m] < \frac{1}{m}$. Since $|X_{m+1}-X_m| = 1$ almost surely, we see that $|E[X_{T_m}]-E[X_{S_m}]| \leq E[T_m -S_m] < \frac{1}{m}$. Since $\sup_{0 \leq m' \leq S_m} X_{m'} \leq k+H_m$, we see by our previous work that $\lim_{m \lar \ff}E[X_{S_m}] = \frac{x_k}{t_\ff}$, and it follows that the same is true for $\lim_{m \lar \ff}E[X_{T_m}]$. This completes the proof of the theorem.

\section{Examples}

{\bf Example 1:} This example appeared in \cite{mebaby}. Let $l_n=\frac{n}{2n+1}, r_n = \frac{n+1}{2n+1}$ for $n \geq 1$. Then $t_n = \frac{1}{n+1}$, so that $t_\ff = 0$. On the other hand, $x_\ff = 1 + \sum_{n=1}^{\ff}t_n = \ff$. We see that the birth-death chain $X_m$ built upon these transition probabilities has an extinction probability of 1, but $E[X_m] \lar \ff$ as $m \lar \ff$.

\vski

{\bf Example 2:} Suppose there is an $M \geq 1$ such that $t_n=r_n=1/2$ for all $n > M$. Then $t_\ff = \frac{l_1 \ldots l_M}{r_1 \ldots r_M}$, and so for a sequence $\{T_m\}_{m=1}^\ff$ satisfying the requirements for our theorem we have

\be
E[X_{T_m}] \lar \frac{1 + \frac{l_1}{r_1} + \ldots + \frac{l_1\ldots l_{k-1}}{r_1\ldots r_{k-1}}}{\frac{l_1 \ldots l_M}{r_1 \ldots r_M}}
\ee

Essentially, if $X_m$ reaches large $n$ it performs simple random walk, and for that reason it is not a surprise that it behaves like a martingale for large stopping times.

\vski

{\bf Example 3:} In order for $t_\ff$ to exist and lie in $(0,\ff)$, it is clearly necessary that $r_n,l_n \lar 1/2$ as $n\lar \ff$, but is not sufficient. Example 1 shows this, as $t_\ff = 0$. Exchanging $r_n$ and $l_n$ in Example 1 will give the case when $t_\ff = \ff$. It is not hard to interlace these series in a way so that $t_\ff$ does not exist. The theory of infinite products(most complex analysis texts contain this, for example see \cite{ahlf}) gives us the standard necessary condition for $t_\ff$ to exist and lie in $(0,\ff)$. The following proposition is obtained.

\begin{proposition} \label{}
Suppose that

\begin{equation} \label{}
\sum_{n=1}^{\ff} \Big| 1- \frac{l_n}{r_n} \Big| < \ff
\end{equation}

Then $\lim_{m \lar \ff}E[X_{T_m}]$ exists and lies in $(0,\ff)$, where $\{T_m\}_{m=1}^\ff$ is any non-decreasing sequence of stopping times satisfying $E[T_m]<\ff$ such that $T_m \lar T_\DD$ almost surely.
\end{proposition}

\section{Concluding remarks}

The reader patient enough to study the proof given here as well as the one given in \cite{mebaby} will no doubt recognize that each is simply the other in disguise. As alluded to in the introduction, the one given here has the advantage of being free of stochastic calculus, which may appeal to the sensibilities of some readers. On the other hand, the other argument shows the local time of Brownian motion to be a very useful tool for bookkeeping, and it is doubtful that the author would have noticed this theorem without studying the Brownian motion model. Regardless, it would be wonderful to have a truly different proof, for instance using electric resistance or perhaps a cleverly constructed change of measure.

\section{Acknowledgements}

The author is grateful for support from Australian Research Council Grant DP0988483.

\bibliographystyle{alpha}
\bibliography{BDChain}

\begin{thebibliography}{Marv1}

\bibitem[Ahl78]{ahlf}
L.V. Ahlfors.
\newblock {\em {Complex Analysis}}.
\newblock McGraw-Hill, 1978.

\bibitem[DS84]{doysne}
P.G. Doyle and J.L. Snell.
\newblock {\em {Random walks and electric networks}}.
\newblock Mathematical Association of America, 1984.

\bibitem[Lan97]{lang}
S.~Lang.
\newblock {\em {Undergraduate analysis}}.
\newblock Springer Verlag, 1997.

\bibitem[Marv1]{mebaby}
G.~Markowsky.
\newblock {Applying Brownian motion to the study of birth-death chains}.
\newblock {\em Statistics and Probability Letters}, to appear,
  arXiv:1103.4208v1.

\bibitem[Nor98]{norris}
J.R. Norris.
\newblock {\em {Markov chains}}.
\newblock Cambridge Univ Pr, 1998.

\bibitem[Sys92]{syski}
R.~Syski.
\newblock {\em {Passage times for Markov chains}}.
\newblock Ios Pr Inc, 1992.

\end{thebibliography}
\end{document}